\pdfoutput=1
\documentclass[11pt,reqno]{amsart}

\usepackage[T1]{fontenc}
\usepackage[utf8]{inputenc}
\usepackage{amsmath,amssymb,amsthm}
\usepackage[letterpaper,margin=0.75in]{geometry}
\usepackage{booktabs}
\usepackage[hidelinks]{hyperref}

\DeclareMathOperator{\vol}{vol}
\DeclareMathOperator{\conv}{conv}
\DeclareMathOperator{\aff}{aff}

\DeclareMathOperator{\interior}{int}

\newcommand{\repobase}{https://github.com/csabamedgyeswork/inscribed-polytopes}
\newcommand{\repocommit}{maxvol-d4-n10-v1}
\newcommand{\repohref}{\href{\repobase}{\texttt{\repobase}}}
\newcommand{\artpath}[1]{\href{\repobase/tree/\repocommit/#1}{\nolinkurl{#1}}}

\newcounter{artifact}
\newenvironment{artifactlist}
  {\sloppy 
   \begin{list}{[R\arabic{artifact}]}
     {\usecounter{artifact}%
      \setlength{\labelwidth}{3.1em}\setlength{\leftmargin}{3.7em}%
      \setlength{\labelsep}{0.6em}\setlength{\itemsep}{0.15ex}%
      \setlength{\parsep}{0pt}}}
  {\end{list}}
\newcommand{\artifact}[3]{\item\label{art:#1}\artpath{#2}, #3}
\newcommand{\artref}[1]{\hyperref[art:#1]{[R\ref*{art:#1}]}}

\newtheorem{theorem}{Theorem}[section]
\newtheorem{proposition}[theorem]{Proposition}
\newtheorem{lemma}[theorem]{Lemma}
\newtheorem{corollary}[theorem]{Corollary}

\theoremstyle{definition}
\newtheorem{definition}[theorem]{Definition}
\newtheorem{remark}[theorem]{Remark}

\newcommand{\citedname}{}
\newtheorem*{citedstatement}{\citedname}
\newenvironment{cited}[2]
  {\renewcommand{\citedname}{#1}\begin{citedstatement}[#2]}
  {\end{citedstatement}}

\title[The inscribed $4$-polytope of maximal volume with ten vertices]
      {The inscribed $4$-polytope of maximal volume\\ with ten vertices}
\author{Csaba Medgyes}
\email{csaba.medgyes.work@gmail.com}

\begin{document}

\begin{abstract}
  We prove that a $4$-polytope with ten vertices inscribed in the unit sphere
  $\mathbb{S}^{3}$ has volume at most $\frac{25}{384}(10 + 2\sqrt{5})$, and
  that, up to orthogonal transformations, equality holds only for the convex
  hull of two regular pentagons inscribed in $\mathbb{S}^{3}$ in orthogonal
  $2$-planes through the origin. As for six to nine vertices, the optimum is
  thus the convex hull of two regular polygons in orthogonal planes. On the
  combinatorial type of the pair, whose boundary complex is the join of two
  $5$-cycles, the sharp bound follows from our work on nine vertices. Since an
  optimal polytope is simplicial, by a result of Horv\'ath and L\'angi, it
  remains to exclude the other $162003$ simplicial types. We do this by
  sum-of-squares certificates checked in exact arithmetic, extended by
  nonnegative multiples of facet polynomials, which are positive on every
  positively oriented realisation of a type. Such a convexity certificate bounds the volume of the
  polytopes of a type without bounding its signed volume on every
  configuration of points on the sphere.
\end{abstract}

\maketitle

\section{Introduction}
\label{sec:intro}

Which convex polytope with $n$ vertices, inscribed in the unit sphere
$\mathbb{S}^{d-1}$ of $\mathbb{R}^{d}$, has the largest volume? Fejes T\'{o}th
raised this question in three dimensions \cite{fejes-toth-1964}, and it is
still open there for most $n$ (see \cite[Section~1]{horvath-langi-2016}).
There, besides the classical cases $n = 4, 6,
12$ (see \cite[Section~1.1]{hoehner-ledford-2026}), Berman and Hanes
determined the optimal polyhedra with at most eight vertices
\cite{berman-hanes-1970}, and Hoehner and Ledford recently completed the case
of nine \cite{hoehner-ledford-2026}, building on a reduction in Wimmer's
thesis \cite{wimmer-1998}. For more vertices only numerical results are
known: Mutoh reported a computer search for up to $30$ vertices
\cite{mutoh-2003}, and Hardin, Sloane and Smith tabulate putatively optimal
configurations for up to $130$ \cite{hardin-sloane-smith}.

In higher dimensions even less is known. Horv\'{a}th and L\'{a}ngi extended
the necessary condition for optimality of Berman and Hanes to every dimension,
and with it showed that optimal polytopes are simplicial
\cite[Lemma~1]{horvath-langi-2016}. They determined the optimal polytopes with
$n = d + 1$ vertices, the regular simplices
\cite[Corollary~1]{horvath-langi-2016}; with $n = d + 2$, the convex
hulls of two regular simplices of dimensions $\lfloor d/2\rfloor$ and
$\lceil d/2\rceil$ in orthogonal subspaces
\cite[Theorem~2]{horvath-langi-2016}; and, for odd $d$, with $n = d + 3$
\cite[Theorem~3]{horvath-langi-2016}. For even $d$ their argument for
$n = d + 3$ excludes cyclic polytopes, and they asked whether an optimal
polytope can be cyclic \cite[Problem~2]{horvath-langi-2016}.

In dimension four these results settle $n = 5$ and $n = 6$. In earlier work
we answered the question of Horv\'{a}th and L\'{a}ngi for $d = 4$, which
settles $n = 7$, and showed that the regular cross-polytope is optimal for
$n = 8$ \cite{medgyes-maxvol-d4}. We then settled $n = 9$
\cite{medgyes-maxvol-d4-n9}. The optima for $n = 6$ to $9$ share a structure:
each is the convex hull of two regular polygons with $m$ and $l$ vertices,
inscribed in $\mathbb{S}^{3}$ in orthogonal $2$-planes through the origin,
with $(m, l) = (3, 3)$, $(3, 4)$, $(4, 4)$ and $(4, 5)$. The method bounds
the volume on each combinatorial type of simplicial polytope by a
sum-of-squares certificate \cite{parrilo-2003}, an algebraic identity checked
in exact rational arithmetic, and treats the type of the optimum separately.
To our knowledge, no case with $d = 4$ and $n \geq 10$ has been settled.

This paper treats $n = 10$, where the pattern suggests a pair of regular
pentagons. Two things change from nine vertices. The number of simplicial
types grows from $1142$ \cite{altshuler-bokowski-steinberg-1980} to $162004$
\cite[Theorem~25(i) and Table~3]{firsching-2017}. And the margin is small:
to exclude the types other than that of the pair, we need a bound less than
$2 \cdot 10^{-4}$ below its volume. At this margin the certificates of
\cite{medgyes-maxvol-d4}, which bound the signed volume on every
configuration of points on the sphere, convex or not, fail on one type, the
cyclic polytope $C_{4}(10)$.

We prove that the pair of regular pentagons is, up to orthogonal
transformations, the unique optimum (Theorem~\ref{thm:main}). On the type of
the pair, the join of two $5$-cycles, the sharp bound follows from
\cite[Corollary~5.2]{medgyes-maxvol-d4-n9}. For $C_{4}(10)$ we add
nonnegative multiples of facet polynomials, which are positive on every
realisation of a type, and obtain a \emph{convexity certificate}
\cite[Section~5]{parrilo-2003}: it bounds only the volume of the polytopes of
the type, which is all the problem requires.

\section{Preliminaries}
\label{sec:prelim}

\subsection{Inscribed polytopes}

\begin{cited}{Notation}{{\cite[Section~1]{horvath-langi-2016}}}
  Throughout the paper, by a polytope we mean a convex polytope. The vertex set
  of a polytope $P$ is denoted by $V(P)$. We denote the family of
  $d$-dimensional polytopes, with $n$ vertices and inscribed in the unit sphere
  $\mathbb{S}^{d-1}$, by $\mathcal{P}_{d}(n)$. We denote $d$-dimensional
  volume by $\vol_{d}$, and set
  $v_{d}(n) = \max\{\vol_{d}(P) : P \in \mathcal{P}_{d}(n)\}$.
\end{cited}

\begin{remark}[{\cite[Remark~2.1]{medgyes-maxvol-d4}}]
  \label{rem:attained}
  For $n \geq d + 1$, $v_{d}(n)$ is attained by some
  $P \in \mathcal{P}_{d}(n)$.
\end{remark}

\begin{definition}[{\cite[Definition~2.1, Section~2.2 and Proposition~2.16]{ziegler-1995}}]
  A \emph{facet} of a $d$-polytope is a face of dimension $d-1$, and a
  \emph{ridge} is a face of dimension $d-2$. A polytope is \emph{simplicial} if
  every facet is a simplex. Two polytopes are \emph{combinatorially
  equivalent} if their face lattices are isomorphic; an equivalence class is a
  \emph{combinatorial type}.
\end{definition}

Equivalently, two polytopes are combinatorially equivalent if their
facet--vertex incidence matrices differ only by permutations of rows and
columns \cite[Section~15.1.1]{henk-richter-gebert-ziegler-2017}. For
simplicial polytopes, which is the form in which we use the definition, this
means that a bijection between the vertex sets carries the vertex sets of the
facets of one onto those of the other.

\begin{definition}[{\cite[Definition~1]{horvath-langi-2016}}]
  \label{def:propertyz}
  Let $P \in \mathcal{P}_{d}(n)$ be a $d$-polytope with
  $V(P) = \{p_{1}, p_{2}, \ldots, p_{n}\}$. If for each $i$, there is an open
  set $U_{i} \subset \mathbb{S}^{d-1}$ such that $p_{i} \in U_{i}$, and for any
  $q \in U_{i}$, we have
  \begin{equation*}
    \vol_{d}\bigl(\conv((V(P) \setminus \{p_{i}\}) \cup \{q\})\bigr)
      \leq \vol_{d}(P),
  \end{equation*}
  then we say that $P$ satisfies Property Z.
\end{definition}

\begin{remark}[{\cite{berman-hanes-1970}, \cite[Remark~2.4]{medgyes-maxvol-d4}}]
  \label{rem:propertyz}
  Every $P \in \mathcal{P}_{d}(n)$ with $\vol_{d}(P) = v_{d}(n)$ satisfies
  Property Z.
\end{remark}

\begin{lemma}[{\cite[Lemma~1]{horvath-langi-2016}}]
  \label{lem:simplicial}
  Every $P \in \mathcal{P}_{d}(n)$ satisfying Property Z is simplicial.
\end{lemma}

\subsection{Signed volume and certificates}

For $v_{1}, \ldots, v_{d} \in \mathbb{R}^{d}$, $\det(v_{1}, \ldots, v_{d})$
denotes the determinant of the matrix with columns $v_{1}, \ldots, v_{d}$.

\begin{definition}[{\cite[Definition~3.1]{medgyes-maxvol-d4}}]
  \label{def:signed-volume}
  An \emph{orientation} of a finite set is an ordering of its elements up to
  even permutations. Let $\mathcal{C} = (n, \mathcal{F})$, where
  $\mathcal{F}$ is a finite set of ordered $d$-tuples
  $F = (F_{1}, \ldots, F_{d})$ of distinct elements of $\{1, \ldots, n\}$,
  called facets. Deleting $F_{a}$ gives a ridge, with the orientation of
  $(F_{1}, \ldots, \widehat{F_{a}}, \ldots, F_{d})$ if $a$ is even and the
  opposite one if $a$ is odd. The orientation of $\mathcal{C}$ is
  \emph{coherent} if every ridge lies in exactly two facets and these induce
  opposite orientations on it. A point
  $p = (p_{1}, \ldots, p_{n}) \in (\mathbb{R}^{d})^{n}$ \emph{realises}
  $\mathcal{C}$ if $P = \conv\{p_{1}, \ldots, p_{n}\}$ has the $p_{i}$ as its
  vertices and its facets are exactly the simplices
  $\conv\{p_{j} : j \in F\}$, $F \in \mathcal{F}$. The \emph{signed volume}
  of $p$ is
  \begin{equation}
    \label{eq:signed-volume}
    D(p) = \sum_{F \in \mathcal{F}} \det(p_{F_{1}}, \ldots, p_{F_{d}}).
  \end{equation}
\end{definition}

We write $j \in F$ if $j$ is one of the entries of $F$, and $p_{F}$ for the
set $\{p_{F_{1}}, \ldots, p_{F_{d}}\}$.

\begin{lemma}[{\cite[Lemma~3.2]{medgyes-maxvol-d4}}]
  \label{lem:volume}
  If the orientation of $\mathcal{C}$ is coherent and $p$ realises
  $\mathcal{C}$, then
  \[
    d!\,\vol_{d}(\conv\{p_{1}, \ldots, p_{n}\}) = |D(p)|.
  \]
\end{lemma}

\begin{definition}[{\cite[Definition~3.4]{medgyes-maxvol-d4}}]
  \label{def:certificate}
  Let $x_{i} = (x_{i,1}, \ldots, x_{i,d})$ be the coordinates of $p_{i}$, and
  $v(x)$ the vector of all monomials of degree at most $2$ in the $x_{i,a}$. A
  \emph{certificate} for $\mathcal{C}$ and $\beta \in \mathbb{Q}$ is a
  symmetric positive definite matrix $Q$ with rational entries together with
  polynomials $h_{1}, \ldots, h_{n}$ with rational coefficients such that
  \begin{equation}
    \label{eq:certificate}
    \beta - D(x) = v(x)^{\mathsf{T}} Q\, v(x)
      + \sum_{i=1}^{n} h_{i}(x)\bigl(|x_{i}|^{2} - 1\bigr)
  \end{equation}
  holds identically as polynomials in the $x_{i,a}$.
\end{definition}

\begin{proposition}[{\cite[Proposition~3.5]{medgyes-maxvol-d4}}]
  \label{prop:certificate}
  Suppose the orientation of $\mathcal{C}$ is coherent and a certificate for
  $\mathcal{C}$ and $\beta$ exists. Then $|D(p)| \leq \beta$ for every
  $p \in (\mathbb{S}^{d-1})^{n}$, and every polytope inscribed in
  $\mathbb{S}^{d-1}$ that realises $\mathcal{C}$ has volume at most
  $\beta/d!$.
\end{proposition}

The first conclusion of Proposition~\ref{prop:certificate} concerns every
configuration on the sphere, convex or not. A certificate therefore cannot
exist for $\beta$ below the maximum of $|D|$ on
$(\mathbb{S}^{d-1})^{n}$, even if that maximum is attained only by
configurations that realise no polytope of the type.
Section~\ref{sec:convexity} removes this restriction.

\subsection{The join of two pentagons}

We use the following results of \cite{medgyes-maxvol-d4-n9} with
$m = l = 5$. For $m, l \geq 3$, $\mathcal{J}_{m,l}$ denotes the join of an
$m$-cycle and an $l$-cycle \cite[Definition~3.1]{medgyes-maxvol-d4-n9}, a
coherently oriented type with $m + l$ vertices
\cite[Lemma~3.3]{medgyes-maxvol-d4-n9}.

\begin{lemma}[{\cite[Lemma~3.2]{medgyes-maxvol-d4-n9}}]
  \label{lem:freesum}
  Let $m, l \geq 3$, let $u_{1}, u_{2}, u_{3}, u_{4}$ be an orthonormal basis of
  $\mathbb{R}^{4}$, and for $i \in \mathbb{Z}_{m}$, $j \in \mathbb{Z}_{l}$ let
  \[
    a_{i} = \cos\tfrac{2\pi i}{m}\,u_{1} + \sin\tfrac{2\pi i}{m}\,u_{2}, \qquad
    b_{j} = \cos\tfrac{2\pi j}{l}\,u_{3} + \sin\tfrac{2\pi j}{l}\,u_{4}.
  \]
  Then $P = \conv\{a_{i}, b_{j}\} \in \mathcal{P}_{4}(m + l)$, and the point
  $p = (a_{1}, \ldots, a_{m}, b_{1}, \ldots, b_{l})$ realises
  $\mathcal{J}_{m,l}$.
\end{lemma}

\begin{corollary}[{\cite[Corollary~5.2]{medgyes-maxvol-d4-n9}}]
  \label{cor:pair}
  Let $m, l \in \{3, 4, 5\}$. If $P \in \mathcal{P}_{4}(m+l)$ is
  combinatorially equivalent to the polytope of Lemma~\ref{lem:freesum}, then
  \[
    \vol_{4}(P) \leq \frac{1}{24}\, m\sin\frac{2\pi}{m}\; l\sin\frac{2\pi}{l},
  \]
  with equality if and only if $P$ is the polytope of Lemma~\ref{lem:freesum}
  for some orthonormal basis $u_{1}, \ldots, u_{4}$, that is, up to an
  orthogonal transformation, the convex hull of a regular $m$-gon and a
  regular $l$-gon inscribed in $\mathbb{S}^{3}$ in orthogonal $2$-planes
  through the origin.
\end{corollary}

For $m = l = 5$ the bound is $V^{*} = \frac{25}{24}\sin^{2}\frac{2\pi}{5}$.
Since $\sin^{2}\frac{2\pi}{5} = \frac{1}{16}(10 + 2\sqrt{5})$,
\begin{equation}
  \label{eq:vstar}
  24\,V^{*} = \frac{25}{16}\bigl(10 + 2\sqrt{5}\bigr), \qquad
  V^{*} = \frac{25}{384}\bigl(10 + 2\sqrt{5}\bigr) = 0.942196\ldots.
\end{equation}

\section{Convexity certificates}
\label{sec:convexity}

Throughout this section $\mathcal{C} = (n, \mathcal{F})$ is as in
Definition~\ref{def:signed-volume}, with $d$-tuples as facets.

\subsection{Facet polynomials}

\begin{definition}
  \label{def:facet-polynomial}
  For $F \in \mathcal{F}$ and $j \in \{1, \ldots, n\}$ with $j \notin F$, the
  \emph{facet polynomial} is
  \begin{equation}
    \label{eq:facet-polynomial}
    A_{F,j}(x) = \det\bigl(x_{F_{1}} - x_{j}, \ldots, x_{F_{d}} - x_{j}\bigr).
  \end{equation}
\end{definition}

$A_{F,j}$ is a polynomial of degree $d$ in the $x_{i,a}$ with integer
coefficients, and $|A_{F,j}(p)|$ is $d!$ times the volume of the simplex
$\conv(p_{F} \cup \{p_{j}\})$. An odd permutation of the entries of $F$
changes its sign, so $A_{F,j}$ depends on the ordering of $F$ given in
$\mathcal{F}$, up to even permutations, and not only on its vertex set.

\begin{lemma}
  \label{lem:facet-sign}
  Suppose the orientation of $\mathcal{C}$ is coherent, $p$ realises
  $\mathcal{C}$, and $D(p) > 0$. Then $A_{F,j}(p) > 0$ for every
  $F \in \mathcal{F}$ and every $j \notin F$.
\end{lemma}

\begin{proof}
  Let $P = \conv\{p_{1}, \ldots, p_{n}\}$, and for $F \in \mathcal{F}$ and
  $y \in \mathbb{R}^{d}$ put
  \[
    a_{F}(y) = \det\bigl(p_{F_{1}} - y, \ldots, p_{F_{d}} - y\bigr),
  \]
  so that $A_{F,j}(p) = a_{F}(p_{j})$. Since $p$ realises $\mathcal{C}$,
  $F \mapsto \conv p_{F}$ maps $\mathcal{F}$ onto the facets of $P$. It is
  injective: two tuples with the same vertex set would share all their
  ridges, and a ridge of one of them would lie in at least three tuples,
  the two and the other tuple at that ridge, against coherence.

  \emph{Step 1: $D(p) = \sum_{F \in \mathcal{F}} a_{F}(y)$ for every
  $y \in \mathbb{R}^{d}$.} Expand $a_{F}(y)$ by multilinearity, writing each
  column $p_{F_{a}} - y$ as the sum of $p_{F_{a}}$ and $-y$. Terms with $-y$
  in two columns vanish, since two of their columns are equal. The term
  without $y$ is $\det(p_{F_{1}}, \ldots, p_{F_{d}})$. The term with $-y$ in
  column $a$ is, after moving $y$ to the first column past $a - 1$ columns,
  \[
    -(-1)^{a-1}\det\bigl(y, p_{F_{1}}, \ldots, \widehat{p_{F_{a}}}, \ldots,
    p_{F_{d}}\bigr)
    = (-1)^{a}\det\bigl(y, p_{F_{1}}, \ldots, \widehat{p_{F_{a}}}, \ldots,
    p_{F_{d}}\bigr).
  \]
  For an ordered $(d-1)$-tuple $R$ put
  $\phi_{y}(R) = \det(y, p_{R_{1}}, \ldots, p_{R_{d-1}})$. It is unchanged
  by even permutations of $R$ and changes sign under odd ones, so it is a
  function of the oriented ridge that changes sign when the orientation is
  reversed. By Definition~\ref{def:signed-volume}, the ridge obtained from
  $F$ by deleting $F_{a}$ carries the orientation of
  $(F_{1}, \ldots, \widehat{F_{a}}, \ldots, F_{d})$ if $a$ is even and the
  opposite one if $a$ is odd, so the term above is $\phi_{y}$ of that
  oriented ridge. By coherence, every ridge lies in exactly two facets,
  which induce opposite orientations on it, so in $\sum_{F} a_{F}(y)$ these
  terms cancel in pairs. What remains is
  $\sum_{F} \det(p_{F_{1}}, \ldots, p_{F_{d}}) = D(p)$.

  \emph{Step 2: $a_{F}(y) > 0$ for all $F \in \mathcal{F}$ and all
  $y \in \interior P$.} Fix $y \in \interior P$ and put
  $S_{F} = \conv(\{y\} \cup p_{F})$, so that $|a_{F}(y)| = d!\,\vol_{d}(S_{F})$.
  The simplices $S_{F}$ cover $P$: for $z \in P \setminus \{y\}$, the ray from
  $y$ through $z$ leaves $P$ at a unique point $w \in \partial P$ with
  $z \in [y, w]$, and $w$ lies in some facet $\conv p_{F}$, so $z \in S_{F}$.
  Their interiors are pairwise disjoint: a point $z$ of the interior of
  $S_{F}$ is $z = (1 - t)y + tu$ with $t \in (0, 1)$ and $u$ in the relative
  interior of $\conv p_{F}$, so $u$ is the point $w$ where the ray from $y$
  through $z$ leaves $P$; if $z$ were also interior to $S_{G}$, then $w$
  would lie in the relative interiors of the facets $\conv p_{F}$ and
  $\conv p_{G}$, and distinct faces of a polytope have disjoint relative
  interiors, so $F = G$. Hence, using Step~1 and Lemma~\ref{lem:volume},
  \[
    \Bigl|\sum_{F} a_{F}(y)\Bigr| = |D(p)| = d!\,\vol_{d}(P)
    = \sum_{F} d!\,\vol_{d}(S_{F}) = \sum_{F} |a_{F}(y)|.
  \]
  Equality in the triangle inequality means that the nonzero $a_{F}(y)$
  all have the same sign. None is zero: $a_{F}(y) = 0$ would put $y$ in the
  affine hull of $p_{F}$, a supporting hyperplane of $P$, which misses
  $\interior P$. Since their sum $D(p)$ is positive, every $a_{F}(y)$ is
  positive.

  \emph{Step 3: conclusion.} Fix $F \in \mathcal{F}$ and $j \notin F$. By the
  expansion in Step~1, $a_{F}$ is an affine function of $y$. It vanishes at
  each $p_{F_{k}}$, since then column $k$ is zero, and it is not constant, by
  Step~2; so its zero set is the hyperplane $H_{F} = \aff p_{F}$. As
  $\conv p_{F}$ is a facet of $P$, $P$ lies in one closed half-space of
  $H_{F}$ and $P \cap H_{F} = \conv p_{F}$. By Step~2 that half-space is
  $\{a_{F} \geq 0\}$, so $a_{F}(p_{j}) \geq 0$. If $a_{F}(p_{j}) = 0$, then
  $p_{j} \in P \cap H_{F} = \conv p_{F}$, so the vertex $p_{j}$ of $P$
  would be a vertex of the facet $\conv p_{F}$, that is, $p_{j} = p_{F_{k}}$
  for some $k$; this is impossible, since $j \notin F$ and the $n$ vertices
  $p_{i}$ of $P$ are distinct. Hence $A_{F,j}(p) = a_{F}(p_{j}) > 0$.
\end{proof}

Only $A_{F,j}(p) \geq 0$ is used below. The hypothesis $D(p) > 0$ fixes the
sign: replacing $p$ by its mirror image changes the signs of $D$ and of
every $A_{F,j}$.

\subsection{Convexity certificates}

\begin{definition}
  \label{def:convexity-certificate}
  Let $v_{1}(x), \ldots, v_{r}(x)$ be vectors of monomials of degree at most
  $2$ in the $x_{i,a}$. A \emph{convexity certificate} for $\mathcal{C}$ and
  $\beta \in \mathbb{Q}$ consists of symmetric positive definite matrices
  $Q_{1}, \ldots, Q_{r}$ with rational entries, polynomials
  $h_{1}, \ldots, h_{n}$ with rational coefficients, and rationals
  $\lambda_{F,j} \geq 0$ for $F \in \mathcal{F}$ and $j \notin F$, such that
  \begin{equation}
    \label{eq:convexity-certificate}
    \beta - D(x)
    = \sum_{b=1}^{r} v_{b}(x)^{\mathsf{T}} Q_{b}\, v_{b}(x)
      + \sum_{F \in \mathcal{F}} \sum_{j \notin F} \lambda_{F,j}\, A_{F,j}(x)
      + \sum_{i=1}^{n} h_{i}(x)\bigl(|x_{i}|^{2} - 1\bigr)
  \end{equation}
  holds identically as polynomials in the $x_{i,a}$.
\end{definition}

With $r = 1$, $v_{1} = v$ and every $\lambda_{F,j} = 0$, a convexity
certificate is a certificate in the sense of Definition~\ref{def:certificate}.
Several blocks $v_{b}$ allow the block-diagonal Gram matrices used in the
computations \cite[Section~3.4]{medgyes-maxvol-d4}; a block may omit
monomials, and two blocks may share one.

\begin{proposition}
  \label{prop:convexity-certificate}
  Suppose the orientation of $\mathcal{C}$ is coherent and a convexity
  certificate for $\mathcal{C}$ and $\beta$ exists. Then every polytope
  inscribed in $\mathbb{S}^{d-1}$ that realises $\mathcal{C}$ has volume at
  most $\beta/d!$.
\end{proposition}

\begin{proof}
  Let $P$ be such a polytope, so that $P = \conv\{p_{1}, \ldots, p_{n}\}$
  for some $p \in (\mathbb{S}^{d-1})^{n}$ that realises $\mathcal{C}$. By
  Lemma~\ref{lem:volume}, $|D(p)| = d!\,\vol_{d}(P) > 0$. If $D(p) < 0$,
  let $\rho$ be a reflection in a hyperplane through the origin, and
  replace $p$ by $\rho p = (\rho p_{1}, \ldots, \rho p_{n})$. Its points
  lie on $\mathbb{S}^{d-1}$, and it realises $\mathcal{C}$, since $\rho$
  maps $P$ onto $\conv\{\rho p_{1}, \ldots, \rho p_{n}\}$ and faces onto
  faces. Moreover $\vol_{d}(\rho P) = \vol_{d}(P)$ and
  $D(\rho p) = \det(\rho)\,D(p) = -D(p)$. So we may assume $D(p) > 0$,
  and then $D(p) = d!\,\vol_{d}(P)$.

  Evaluate \eqref{eq:convexity-certificate} at $p$. The last sum vanishes
  since $|p_{i}| = 1$. Each $v_{b}(p)^{\mathsf{T}} Q_{b}\, v_{b}(p)$ is
  nonnegative since $Q_{b}$ is positive definite. Each
  $\lambda_{F,j} A_{F,j}(p)$ is nonnegative since $\lambda_{F,j} \geq 0$ and,
  by Lemma~\ref{lem:facet-sign}, $A_{F,j}(p) > 0$. Hence
  $\beta - d!\,\vol_{d}(P) \geq 0$, that is,
  $\vol_{d}(P) \leq \beta/d!$.
\end{proof}

\begin{remark}
  \label{rem:what-it-bounds}
  Let $K_{\mathcal{C}}$ be the set of $p \in (\mathbb{S}^{d-1})^{n}$ with
  $A_{F,j}(p) \geq 0$ for all $F \in \mathcal{F}$ and $j \notin F$.
  Identity \eqref{eq:convexity-certificate} is a sum-of-squares certificate
  for $D \leq \beta$ on the semialgebraic set $K_{\mathcal{C}}$
  \cite[Section~5]{parrilo-2003}. It says nothing about $D$ outside
  $K_{\mathcal{C}}$, nor about $-D$. The proof above needs neither, since
  every realisation has a mirror image with $D > 0$, which lies in
  $K_{\mathcal{C}}$ by Lemma~\ref{lem:facet-sign}. So, unlike a certificate
  of Definition~\ref{def:certificate}, a convexity certificate can exist for
  $\beta$ below the maximum of $|D|$ on $(\mathbb{S}^{d-1})^{n}$, but not
  below the maximum of $D$ on $K_{\mathcal{C}}$. That maximum is at least
  $d!$ times the largest volume of an inscribed realisation, and may be
  larger, since $K_{\mathcal{C}}$ also contains configurations that realise
  no polytope, for instance ones in which two points coincide.
\end{remark}

\subsection{Verification}
\label{sec:verification}

A convexity certificate is finite rational data: the ordered facets of
$\mathcal{C}$, $\beta$, the blocks $v_{b}$ and entries of the $Q_{b}$, the
coefficients of the $h_{i}$, and the triples $(F, j, \lambda_{F,j})$ with
$\lambda_{F,j} \neq 0$. The verifier checks the hypotheses of
Proposition~\ref{prop:convexity-certificate} in exact rational arithmetic:
\begin{enumerate}
  \item \emph{Multipliers.} For every stored triple, $\lambda_{F,j} \geq 0$,
    $F$ is one of the stored facets with its stored ordering, and
    $j \in \{1, \ldots, n\}$ is not an entry of $F$. The ordering matters,
    since reversing it changes the sign of $A_{F,j}$.
  \item \emph{Identity.} $D$ is rebuilt from the stored facets by
    \eqref{eq:signed-volume} and each $A_{F,j}$ by
    \eqref{eq:facet-polynomial}, both sides of
    \eqref{eq:convexity-certificate} are expanded, and every coefficient
    must agree.
  \item \emph{Positive definiteness.} Each $Q_{b}$ is symmetric, and either
    $M = R^{\mathsf{T}} Q_{b} R$ is strictly diagonally dominant with
    positive diagonal for some upper triangular $R$ with nonzero diagonal
    entries, or $Q_{b}$ has positive leading principal minors, computed by
    the fraction-free elimination of Bareiss \cite{bareiss-1968}. In the
    first case $M$ is positive definite
    \cite[Theorem~6.1.10(c)]{horn-johnson-2013}, and hence so is $Q_{b}$, since
    $R$ is invertible; in the second, Sylvester's criterion
    \cite[Theorem~7.2.5]{horn-johnson-2013} applies. The verifier takes for
    $R$ a rounded floating-point inverse Cholesky factor of $Q_{b}$, with
    integer entries, and computes $M$ exactly; the rounding only chooses $R$.
  \item \emph{Orientation and type.} Every ridge lies in exactly two stored
    facets, with opposite induced orientations, so the orientation is
    coherent; and the vertex sets of the stored facets are exactly those of
    the combinatorial type the certificate is for.
\end{enumerate}
Checks 1 to 3 show that the data is a convexity certificate, and check~4
that $\mathcal{C}$ is coherent and is the intended type. Check~4 needs the
facets of the type to be exact: in Section~\ref{sec:reduction} they come
from polytopes with rational vertices, by exact supporting-hyperplane tests
as in \cite[Lemma~6.1]{medgyes-maxvol-d4-n9}.

\section{The other types and the main theorem}
\label{sec:reduction}

\begin{theorem}[{\cite[Theorem~25(i) and Table~3]{firsching-2017}}]
  \label{thm:count}
  There are exactly $162004$ combinatorial types of simplicial $4$-polytopes
  with $10$ vertices.
\end{theorem}

Firsching's data \cite{firsching-2017} give a polytope with rational
vertices for each line of a file with $162004$ lines \artref{data}; we label
them \texttt{firsching-10-0} to \texttt{firsching-10-162003}.

\begin{lemma}
  \label{lem:data}
  The polytopes \texttt{firsching-10-0} to \texttt{firsching-10-162003}
  represent every combinatorial type of simplicial $4$-polytope with ten
  vertices, each exactly once, and the type of \texttt{firsching-10-11} is
  that of the polytope of Lemma~\ref{lem:freesum} with $m = l = 5$.
\end{lemma}

\begin{proof}
  We check in exact arithmetic, as for nine vertices
  \cite[Lemma~6.1]{medgyes-maxvol-d4-n9}, that these $162004$ polytopes are
  simplicial and pairwise not combinatorially equivalent \artref{catalogue},
  \artref{testten}: the facets of each are the $4$-subsets of its vertices
  with every other vertex strictly on one side of their hyperplane, every
  $3$-subset of a facet lies in exactly two facets, and the $162004$ facet
  systems have pairwise distinct canonical forms under relabelling. By
  Theorem~\ref{thm:count} they represent every type. The facets of
  \texttt{firsching-10-11} are exactly the sets
  $\{a_{i}, a_{i+1}, b_{j}, b_{j+1}\}$, indices modulo $5$, for the
  $5$-cycles $(a_{0}, \ldots, a_{4}) = (0, 2, 7, 5, 6)$ and
  $(b_{0}, \ldots, b_{4}) = (1, 4, 3, 8, 9)$ of its zero-based vertex labels
  \artref{testten}. These are the facets of $\mathcal{J}_{5,5}$, and by
  Lemma~\ref{lem:freesum} so are those of the polytope there, so the two are
  combinatorially equivalent.
\end{proof}

The polytopes need not be inscribed: a convexity certificate is a polynomial
identity and needs no realisation, and for a type without an inscribed
realisation the conclusion of Proposition~\ref{prop:convexity-certificate}
holds vacuously.

\begin{proposition}
  \label{prop:rivals}
  Every simplicial $P \in \mathcal{P}_{4}(10)$ that is not combinatorially
  equivalent to \texttt{firsching-10-11} has
  $\vol_{4}(P) \leq 471/500 < V^{*}$, with $V^{*}$ as in \eqref{eq:vstar}.
\end{proposition}

\begin{proof}
  By Lemma~\ref{lem:data}, $P$ is combinatorially equivalent to
  \texttt{firsching-10-$k$} for some $k \neq 11$. For each of these $162003$
  values of $k$, the catalogue \artref{catalogue} gives a coherently oriented
  $\mathcal{C}_{k} = (10, \mathcal{F}_{k})$ whose facet vertex sets are those
  of \texttt{firsching-10-$k$}, and there is a convexity certificate for
  $\mathcal{C}_{k}$ and $\beta = 2826/125$ that passes the exact checks of
  Section~\ref{sec:verification} \artref{verifier}, \artref{testten}
  (Remark~\ref{rem:computations}).

  Since $P$ is simplicial, the combinatorial equivalence is a bijection
  between the vertex sets that maps the facet vertex sets of
  \texttt{firsching-10-$k$} onto those of $P$. Numbering the vertices of $P$
  through it as $p = (p_{1}, \ldots, p_{10})$, $p$ realises
  $\mathcal{C}_{k}$. Proposition~\ref{prop:convexity-certificate} gives
  $\vol_{4}(P) \leq \beta/4! = 471/500$, and $471/500 < V^{*}$ by
  \eqref{eq:vstar}.
\end{proof}

\begin{remark}
  \label{rem:cyclic}
  Only one type needs the facet polynomials. For $162002$ of the $162003$
  values of $k$, the certificate found has every $\lambda_{F,j} = 0$, so it
  is a certificate in the sense of Definition~\ref{def:certificate}
  \artref{ledger}. The exception is \texttt{firsching-10-157156}, the cyclic
  polytope $C_{4}(10)$: its $35$ facet vertex sets are, after relabelling,
  the $4$-subsets of $\{1, \ldots, 10\}$ that satisfy Gale's evenness
  condition \cite[Theorem~0.7]{ziegler-1995}. For this type no certificate
  of Definition~\ref{def:certificate} with $\beta = 2826/125$ exists. Place
  the vertices with zero-based labels $0, 6, 2, 4, 7$ and $1, 9, 3, 5, 8$, in
  this order, at the consecutive vertices of two regular pentagons inscribed
  in $\mathbb{S}^{3}$ in orthogonal $2$-planes through the origin. The convex
  hull is the pair of pentagons, so this point does not realise $C_{4}(10)$,
  and an exact computation gives $|D| = 24 V^{*} > \beta$ there; hence
  Proposition~\ref{prop:certificate} excludes such a certificate. Numerically, $24 V^{*}$ is also the maximum of
  $|D|$ on $(\mathbb{S}^{3})^{10}$ for this type. The convexity certificate
  for it is stored \artref{cyclic}.
\end{remark}

\begin{theorem}
  \label{thm:main}
  \[
    v_{4}(10) = \frac{25}{24}\sin^{2}\frac{2\pi}{5}
      = \frac{25}{384}\bigl(10 + 2\sqrt{5}\bigr) = 0.942196\ldots,
  \]
  and every $P \in \mathcal{P}_{4}(10)$ with $\vol_{4}(P) = v_{4}(10)$ is, up to
  an orthogonal transformation, the convex hull of two regular pentagons
  inscribed in $\mathbb{S}^{3}$ in orthogonal $2$-planes through the origin.
\end{theorem}

\begin{proof}
  By Lemma~\ref{lem:freesum} and Corollary~\ref{cor:pair} with
  $m = l = 5$, the pair of pentagons lies in $\mathcal{P}_{4}(10)$ and has
  volume $V^{*}$, so $v_{4}(10) \geq V^{*}$.

  Let $P \in \mathcal{P}_{4}(10)$ with $\vol_{4}(P) = v_{4}(10)$, which exists
  by Remark~\ref{rem:attained}. By Remark~\ref{rem:propertyz} and
  Lemma~\ref{lem:simplicial}, $P$ is simplicial. Since
  $\vol_{4}(P) \geq V^{*}$, Proposition~\ref{prop:rivals} shows that $P$ is
  combinatorially equivalent to \texttt{firsching-10-11}, and so, by
  Lemma~\ref{lem:data}, to the polytope of Lemma~\ref{lem:freesum} with
  $m = l = 5$. Corollary~\ref{cor:pair} gives
  $v_{4}(10) = \vol_{4}(P) \leq V^{*}$, so $v_{4}(10) = V^{*}$, and its
  equality case identifies $P$.
\end{proof}

\begin{remark}
  \label{rem:computations}
  The certificates are found by semidefinite programming and rounded to
  rationals as in \cite[Section~3.4]{medgyes-maxvol-d4}; the proof uses only
  the exact checks of Section~\ref{sec:verification}. The certificate for
  $C_{4}(10)$ is stored \artref{cyclic}. The other $162002$ are not: the
  ledger \artref{ledger} records the bound and the kind of certificate for
  each type, and running the search again on the exact facets
  \artref{catalogue} rebuilds one. Since floating-point solvers are not
  reproducible bit for bit, a rebuilt certificate may differ from the
  original, but any certificate that passes the exact checks proves the
  bound for its type.
\end{remark}

\section*{AI assistance}

An AI assistant was used during this work, including for mathematical
discussion, drafting and revising the text, and checking the cited statements
against their sources. The proofs do not depend on the assistant: they rest on
the cited results, on the arguments given here, and on the certificates, which
the exact verifier checks.

\section*{Code and data}
\label{sec:artifacts}

{\footnotesize
Everything reported here is available at \repohref, in the version tagged
\href{\repobase/tree/\repocommit}{\texttt{\repocommit}}. The entries below are cited
as [R$k$]; they are artifacts rather than literature and are numbered separately
from the bibliography.

\begin{artifactlist}
  \item\label{art:data}\nolinkurl{data/polytopalsimplicialspheres10_4.txt},
    Firsching's realisations of the $162004$ simplicial types with ten
    vertices, from \url{https://firsching.ch/polytopes}, fetched by
    \artpath{scripts/setup\_data.py}; its SHA-256 is recorded in
    \artpath{results/certificates/d4-n10/manifest.json}
  \artifact{catalogue}{results/certificates/d4-n10/types}{the exact oriented
    facets of the $162004$ polytopes, with
    \artpath{results/certificates/d4-n10/audit.json}, built by
    \artpath{src/inscribed/elimination/catalogue.py}}
  \artifact{ledger}{results/certificates/d4-n10/ledger}{one entry per type,
    with \artpath{results/certificates/d4-n10/manifest.json}, written and
    rebuilt by \artpath{scripts/eliminate\_convex.py}}
  \artifact{cyclic}{results/certificates/d4-n10/firsching-10-157156.convexity.json.gz}{the
    convexity certificate for $C_{4}(10)$ of Remark~\ref{rem:cyclic}}
  \artifact{verifier}{src/inscribed/certify/convexity.py}{the exact verifier of
    Section~\ref{sec:verification}, extending
    \artpath{src/inscribed/certify/verify.py} of
    \cite[Section~3.3]{medgyes-maxvol-d4}, with
    \artpath{src/inscribed/combinatorics/complex.py} checking coherence}
  \artifact{testten}{tests/papers/maxvol\_d4\_n10/test\_ten.py}{tests of the
    claims about Firsching's data and the certificates; the checks over all
    $162004$ types run with \texttt{poetry run pytest -m slow}}
\end{artifactlist}
}

\renewcommand{\bysame}{Csaba Medgyes}
\bibliographystyle{amsplain}
\bibliography{references}

@article{altshuler-bokowski-steinberg-1980,
  author  = {Altshuler, Amos and Bokowski, J\"{u}rgen and Steinberg, Leon},
  title   = {The classification of simplicial $3$-spheres with nine vertices
             into polytopes and nonpolytopes},
  journal = {Discrete Mathematics},
  volume  = {31},
  number  = {2},
  pages   = {115--124},
  year    = {1980},
  doi     = {10.1016/0012-365X(80)90029-1},
}

@article{berman-hanes-1970,
  author  = {Berman, Joel D. and Hanes, Kit},
  title   = {Volumes of polyhedra inscribed in the unit sphere in {$E^3$}},
  journal = {Mathematische Annalen},
  volume  = {188},
  pages   = {78--84},
  year    = {1970},
  doi     = {10.1007/BF01435416},
}

@book{fejes-toth-1964,
  author    = {Fejes T\'{o}th, L\'{a}szl\'{o}},
  title     = {Regular Figures},
  publisher = {Pergamon Press},
  address   = {Oxford},
  year      = {1964},
}

@article{firsching-2017,
  author  = {Firsching, Moritz},
  title   = {Realizability and inscribability for simplicial polytopes via
             nonlinear optimization},
  journal = {Mathematical Programming},
  volume  = {166},
  pages   = {273--295},
  year    = {2017},
  doi     = {10.1007/s10107-017-1120-0},
  eprint  = {1508.02531},
  archivePrefix = {arXiv},
}

@misc{hardin-sloane-smith,
  author       = {Hardin, R. H. and Sloane, N. J. A. and Smith, W. D.},
  title        = {Maximal volume spherical codes},
  howpublished = {\url{https://neilsloane.com/maxvolumes/}},
  url          = {https://neilsloane.com/maxvolumes/},
}

@incollection{henk-richter-gebert-ziegler-2017,
  author    = {Henk, Martin and Richter-Gebert, J\"{u}rgen and Ziegler,
               G\"{u}nter M.},
  title     = {Basic properties of convex polytopes},
  booktitle = {Handbook of Discrete and Computational Geometry},
  edition   = {Third},
  editor    = {Goodman, Jacob E. and O'Rourke, Joseph and T\'{o}th, Csaba D.},
  chapter   = {15},
  pages     = {383--413},
  publisher = {CRC Press},
  address   = {Boca Raton, FL},
  year      = {2017},
}

@misc{hoehner-ledford-2026,
  author = {Hoehner, Steven and Ledford, Jeff},
  title  = {The maximum volume polytope with nine vertices inscribed in the
            sphere},
  year   = {2026},
  howpublished = {arXiv:2608.04392v2 [math.MG]},
  eprint = {2608.04392},
  archivePrefix = {arXiv},
  primaryClass  = {math.MG},
}

@book{horn-johnson-2013,
  author    = {Horn, Roger A. and Johnson, Charles R.},
  title     = {Matrix Analysis},
  edition   = {Second},
  publisher = {Cambridge University Press},
  address   = {Cambridge},
  year      = {2013},
}

@article{horvath-langi-2016,
  author  = {Horv\'{a}th, \'{A}kos G. and L\'{a}ngi, Zsolt},
  title   = {Maximum volume polytopes inscribed in the unit sphere},
  journal = {Monatshefte f\"{u}r Mathematik},
  volume  = {181},
  number  = {2},
  pages   = {341--354},
  year    = {2016},
  doi     = {10.1007/s00605-016-0949-2},
  eprint  = {1402.6496},
  archivePrefix = {arXiv},
}

@misc{medgyes-maxvol-d4,
  author = {Medgyes, Csaba},
  title  = {Inscribed $4$-polytopes of maximal volume with seven and eight
            vertices},
  year   = {2026},
  howpublished = {arXiv:2610.26886 [math.MG]},
  eprint = {2610.26886},
  archivePrefix = {arXiv},
  primaryClass  = {math.MG},
}

@misc{medgyes-maxvol-d4-n9,
  author = {Medgyes, Csaba},
  title  = {The inscribed $4$-polytope of maximal volume with nine vertices},
  year   = {2026},
  howpublished = {arXiv:2611.10078 [math.MG]},
  eprint = {2611.10078},
  archivePrefix = {arXiv},
  primaryClass  = {math.MG},
}

@phdthesis{wimmer-1998,
  author = {Wimmer, Lienhard},
  title  = {Das maximale {V}olumen der konvexen {H}\"{u}lle von {P}unkten auf
            der {E}inheitskugel},
  school = {Universit\"{a}t Salzburg},
  year   = {1998},
}

@incollection{mutoh-2003,
  author    = {Mutoh, Nobuaki},
  title     = {The polyhedra of maximal volume inscribed in the unit sphere and
               of minimal volume circumscribed about the unit sphere},
  booktitle = {Discrete and Computational Geometry (JCDCG 2002)},
  editor    = {Akiyama, Jin and Kano, Mikio},
  series    = {Lecture Notes in Computer Science},
  volume    = {2866},
  pages     = {204--214},
  publisher = {Springer},
  address   = {Berlin, Heidelberg},
  year      = {2003},
  doi       = {10.1007/978-3-540-44400-8_22},
}

@book{ziegler-1995,
  author    = {Ziegler, G\"{u}nter M.},
  title     = {Lectures on Polytopes},
  series    = {Graduate Texts in Mathematics},
  volume    = {152},
  publisher = {Springer},
  address   = {New York},
  year      = {1995},
  doi       = {10.1007/978-1-4613-8431-1},
}

@article{parrilo-2003,
  author  = {Parrilo, Pablo A.},
  title   = {Semidefinite programming relaxations for semialgebraic problems},
  journal = {Mathematical Programming, Series B},
  volume  = {96},
  number  = {2},
  pages   = {293--320},
  year    = {2003},
  doi     = {10.1007/s10107-003-0387-5},
}

@article{bareiss-1968,
  author  = {Bareiss, Erwin H.},
  title   = {Sylvester's identity and multistep integer-preserving {G}aussian
             elimination},
  journal = {Mathematics of Computation},
  volume  = {22},
  number  = {103},
  pages   = {565--578},
  year    = {1968},
  doi     = {10.1090/S0025-5718-1968-0226829-0},
}

\end{document}